\documentclass[a4paper,notitlepage]{article}%
\usepackage{amssymb}
\usepackage{graphicx}
\usepackage{amsmath}
\usepackage{amsthm}
\usepackage{amsfonts}%
\providecommand{\U}[1]{\protect\rule{.1in}{.1in}}
\theoremstyle{plain}
\newtheorem{theorem}{Theorem}[section]

\newtheorem{corollary}[theorem]{Corollary}

\newtheorem{lemma}[theorem]{Lemma}

\theoremstyle{definition}

\newtheorem{remark}[theorem]{Remark}
\numberwithin{equation}{section}
\numberwithin{theorem}{section}
\ifx\pdfoutput\relax\let\pdfoutput=\undefined\fi
\newcount\msipdfoutput
\ifx\pdfoutput\undefined\else
\ifcase\pdfoutput\else
\ifx\paperwidth\undefined\else
\ifdim\paperheight=0pt\relax\else\pdfpageheight\paperheight\fi
\ifdim\paperwidth=0pt\relax\else\pdfpagewidth\paperwidth\fi
\fi\fi\fi
\begin{document}

\title{Positive solutions for nonlinear problems involving the one-dimensional $\phi
$-Laplacian\thanks{2000 \textit{Mathematics Subject Clasification}. 34B15;
34B18; 35J25.} 
\thanks{Partially supported by Secyt-UNC 30720150100019CB.} }
\author{U. Kaufmann, L. Milne\thanks{\textit{E-mail addresses. }%
kaufmann@mate.uncor.edu (U. Kaufmann, Corresponding Author),
milne@mate.uncor.edu (L. Milne).}
\and \noindent\\{\small FaMAF, Universidad Nacional de C\'{o}rdoba, (5000) C\'{o}rdoba,
Argentina}}
\maketitle

\begin{abstract}
Let $\Omega:=\left(  a,b\right)  \subset\mathbb{R}$, $m\in L^{1}\left(
\Omega\right)  $ and $\lambda>0$ be a real parameter. Let $\mathcal{L}$ be the
differential operator given by $\mathcal{L}u:=-\phi\left(  u^{\prime}\right)
^{\prime}+r\left(  x\right)  \phi\left(  u\right)  $, where $\phi
:\mathbb{R\rightarrow R}$ is an odd increasing homeomorphism and $0\leq r\in
L^{1}\left(  \Omega\right)  $. We study the existence of \textit{positive}
solutions for problems of the form%
\[
\left\{
\begin{array}
[c]{ll}%
\mathcal{L}u=\lambda m\left(  x\right)  f\left(  u\right)   & \text{in }%
\Omega,\\
u=0 & \text{on }\partial\Omega,
\end{array}
\right.
\]
where $f:\left[  0,\infty\right)  \rightarrow\left[  0,\infty\right)  $ is a
continuous function which is, roughly speaking, sublinear with respect to
$\phi$. Our approach combines the sub and supersolution method with some
estimates on related nonlinear problems. We point out that our results are new
even in the cases $r\equiv0$ and/or $m\geq0$.

\end{abstract}

\section{Introduction}

Let $\Omega:=\left(  a,b\right)  \subset\mathbb{R}$, $m\in L^{1}\left(
\Omega\right)  $ and $\lambda>0$ be a real parameter. Let us consider problems
of the form%
\begin{equation}
\left\{
\begin{array}
[c]{ll}%
-\phi\left(  u^{\prime}\right)  ^{\prime}=\lambda m\left(  x\right)  f\left(
u\right)  & \text{in }\Omega,\\
u=0 & \text{on }\partial\Omega,
\end{array}
\right.  \label{intro}%
\end{equation}
where $\phi:\mathbb{R\rightarrow R}$ is an odd increasing homeomorphism and
$f:\left[  0,\infty\right)  \rightarrow\left[  0,\infty\right)  $ is a
continuous function. The existence of \textit{positive} solutions for problems
as (\ref{intro}) involving the so-called $\phi$-Laplacian have been widely
studied in the literature (see e.g. \cite{ba, bai, be, bee, ga, ka, lee, luo,
ya} and the references therein) and appear in diverse applications such as
reaction-diffusion systems, nonlinear elasticity, glaciology, population
biology, combustion theory, and non-Newtonian fluids, see for instance
\cite{diaz, dra, glo, o}. We mention also that these kind of problems arise
naturally in the study of radial solutions for nonlinear equations in annular
domains (see e.g. \cite{wang} and its references).

When $\phi\left(  x\right)  =\left\vert x\right\vert ^{p-2}x$ and $f\left(
x\right)  =x^{q}$ with $1<p<\infty$ and $0<q<p-1$, the existence of positive
solutions for (\ref{intro}) was considered in \cite{bams}, even for
sign-changing weights (see also \cite{drabek, or} for the analogous
$N$-dimensional problem). We note, however, that for the computations in
\cite{bams} it was crucial the \textit{homogeneity} of both $\phi$ and $f$,
which of course is no longer true here.

Let us now introduce the following assumptions on $m$ and $\phi$:

\begin{enumerate}
\item[$(M)$] $m\in\mathcal{C}(\overline{\Omega})$ with $m\geq0$ in $\Omega$
and $m\not \equiv 0$ on any subinterval of $\Omega,$

\item[$\left(  M^{\prime}\right)  $] $m\in\mathcal{C}(\overline{\Omega})$ with
$\min_{\overline{\Omega}}m>0,$

\item[$\left(  \Phi\right)  $] There exist increasing homeomorphisms $\psi
_{1},\psi_{2}:\left[  0,\infty\right)  \rightarrow\left[  0,\infty\right)  $
such that $\psi_{1}\left(  t\right)  \phi\left(  x\right)  \leq\phi\left(
tx\right)  \leq\psi_{2}\left(  t\right)  \phi\left(  x\right)  $ for all
$t,x>0,$

\item[$\left(  \Phi^{\prime}\right)  $] There exist $p,q\in(0,\infty)$ such
that $t^{q}\phi\left(  x\right)  \leq\phi\left(  tx\right)  \leq t^{p}%
\phi\left(  x\right)  $ for $t\in\left[  0,1\right]  $ and all $x>0.$
\end{enumerate}

\noindent Under some standard growth conditions on $f$ (which allow both
sublinear and superlinear nonlinearities) and assuming $(M)$ and $(\Phi)$, it
was proved that (\ref{intro}) possesses a positive solution for all
$\lambda>0$ (see \cite[Theorem 1.1]{wang2}), and recently in \cite[Theorem
2]{xu} the authors extended this result to certain $m\in L_{loc}^{1}\left(
\Omega\right)  $ and not requiring that $\psi_{2}\left(  0\right)  =0$. We
point out that these hypothesis impose, in particular, rather strong
restrictions on
\[
l\left(  t\right)  :=\underline{\lim}_{x\rightarrow\infty}\phi\left(
tx\right)  /\phi\left(  x\right)  \quad\text{and}\quad L\left(  t\right)
:=\overline{\lim}_{x\rightarrow\infty}\phi\left(  tx\right)  /\phi\left(
x\right)  .
\]
Indeed, the existence of $\psi_{1}$ as above implies that $l\left(  t\right)
>0$ for all $t\in\left(  0,1\right)  $ and $\lim_{t\rightarrow\infty}l\left(
t\right)  =\infty$, while the existence of $\psi_{2}$ entails that $L\left(
t\right)  <\infty$ for all $t>1$. Let us note that the first and third of
these conditions are not satisfied for instance by exponential-like
nonlinearities, and the remaining one does not hold for example for
logarithmic-like functions.

On the other side, a similar result was established in \cite[Corollary
3.4]{ben} assuming $(M^{\prime})$ and $(\Phi^{\prime})$. We observe that the
first inequality in $(\Phi^{\prime})$ also implies that $l\left(  t\right)
>0$ for all $t\in\left(  0,1\right)  $, while the second one requires that
$\lim_{t\rightarrow0^{+}}L\left(  t\right)  =0$ (and this does not occur, for
instance, with logarithmic-like nonlinearities). Let us add that in all these
works the main tool utilized was some kind of Krasnoselskii-type fixed point
theorem in cones.

Following a different approach, in Theorem \ref{nuevo} below we shall improve
substantially the aforementioned results in the sublinear case, under much
weaker conditions on both $\phi$ and $m$. In fact, regarding the assumptions
on $m\in L^{1}\left(  \Omega\right)  $, we shall only require that $0\leq
m\not \equiv 0$ in $\Omega$. Furthermore, we shall see that the solutions
$u_{\lambda}\rightarrow0$ in $\mathcal{C}^{1}(\overline{\Omega})$ as
$\lambda\rightarrow0^{+}$. In order to derive our theorems, we shall rely on
the well-known sub and supersolution method, combined with upper and lower
estimates on some related nonlinear problems.

Also, under some additional hypothesis on $\phi$ and $m$, we shall prove in
Theorem \ref{r} similar results for the differential operator%
\begin{equation}
\mathcal{L}u:=-\phi\left(  u^{\prime}\right)  ^{\prime}+r\left(  x\right)
\phi\left(  u\right)  ,\label{L}%
\end{equation}
where $0\leq r\in L^{1}\left(  \Omega\right)  $. Moreover, as a consequence of
Theorems \ref{nuevo} and \ref{r}, we shall deduce the existence of
(nontrivial) nonnegative solutions for \textit{sign-changing} weights $m$, see
Corollary \ref{cori}.

The rest of the article is organized as follows. In the next section we
collect some auxiliary results, while in Section 3 we shall state and prove
our main theorems. Finally, at the end of the paper we present several
examples illustrating our conditions and their relations with the ones already
mentioned (see also Remarks \ref{tre} and \ref{lapla}).

\section{Preliminaries}

Let $\phi:\mathbb{R\rightarrow R}$ be an odd increasing homeomorphism and
$h\in L^{1}(\Omega)$. We start compiling some necessary facts about the
problem%
\begin{equation}
\left\{
\begin{array}
[c]{ll}%
-\phi\left(  v^{\prime}\right)  ^{\prime}=h\left(  x\right)  & \text{in
}\Omega,\\
v=0 & \text{on }\partial\Omega.
\end{array}
\right.  \label{g}%
\end{equation}

\begin{remark}
\label{uno}For every $h\in L^{1}(\Omega)$, (\ref{g}) admits a unique solution
$v\in\mathcal{C}^{1}(\overline{\Omega})$ such that $\phi\left(  v^{\prime
}\right)  $ is absolutely continuous and that the equation holds pointwise
$a.e.$ $x\in\Omega$. In fact, one can see that
\begin{equation}
v\left(  x\right)  =\int_{a}^{x}\phi^{-1}\left(  c_{h}-\int_{a}^{y}h\left(
t\right)  dt\right)  dy, \label{v1}%
\end{equation}
where $c_{h}$ is the unique constant such that $v\left(  b\right)  =0$.
Furthermore, the solution operator $\mathcal{S}_{\phi}:L^{1}(\Omega
)\rightarrow\mathcal{C}^{1}(\overline{\Omega})$ is continuous (see e.g.
\cite[Lemma 2.1]{dang}).
\end{remark}

The following lemma shows that $\mathcal{S}_{\phi}$ is a nondecreasing
operator. Although this result should probably be well-known, we have not been
able to find a proof in the literature.

\begin{lemma}
\label{node}Let $h_{1},h_{2}\in L^{1}\left(  \Omega\right)  $ with $h_{1}\leq
h_{2}$ $a.e.$ $x\in\Omega$. Then $\mathcal{S}_{\phi}(h_{1})\leq\mathcal{S}%
_{\phi}(h_{2})$ in $\overline{\Omega}$.
\end{lemma}

\textit{Proof}. Let $v_{i}:=\mathcal{S}_{\phi}(h_{i})$, $i=1,2$, and suppose
by contradiction that $\mathcal{O}:=\left\{  x\in\Omega:v_{1}>v_{2}\right\}
\not =\emptyset$. Let $\mathcal{O}_{c}$ be a connected component of
$\mathcal{O}$. Note that, either by the continuity of $v_{1}$ and $v_{2}$ or
by the boundary condition in (\ref{g}), $v_{1}=v_{2}$ on $\partial
\mathcal{O}_{c}$. Taking into account this, multiplying (\ref{g}) (with
$h_{1}$ in place of $h$) by $v_{1}-v_{2}$ and integrating by parts we get%
\[
\int_{\mathcal{O}_{c}}\phi\left(  v_{1}^{\prime}\right)  \left(  v_{1}%
-v_{2}\right)  ^{\prime}=\int_{\mathcal{O}_{c}}h_{1}\left(  x\right)  \left(
v_{1}-v_{2}\right)  .
\]
Since we can argue in the same way with the equation involving $h_{2}$ and
$h_{1}\leq h_{2}$ in $\Omega$, recalling that $\phi$ is increasing we infer
that%
\[
0\leq\int_{\mathcal{O}_{c}}\left(  \phi\left(  v_{1}^{\prime}\right)
-\phi\left(  v_{2}^{\prime}\right)  \right)  \left(  v_{1}^{\prime}%
-v_{2}^{\prime}\right)  =\int_{\mathcal{O}_{c}}\left(  h_{1}\left(  x\right)
-h_{2}\left(  x\right)  \right)  \left(  v_{1}-v_{2}\right)  \leq0
\]
and thus $v_{1}^{\prime}=v_{2}^{\prime}$ in $\mathcal{O}_{c}$. Furthermore,
$v_{1}=v_{2}$ in $\mathcal{O}_{c}$ because $v_{1}=v_{2}$ on $\partial
\mathcal{O}_{c}$. Contradiction. $\blacksquare$

\quad

For $h\in L^{1}(\Omega)$ with $0\leq h\not \equiv 0$ we define%
\begin{align*}
\mathcal{A}_{h}  &  :=\left\{  x\in\Omega:h\left(  y\right)  =0\text{
}a.e.\text{ }y\in\left(  a,x\right)  \right\}  ,\\
\mathcal{B}_{h}  &  :=\left\{  x\in\Omega:h\left(  y\right)  =0\text{
}a.e.\text{ }y\in\left(  x,b\right)  \right\}  ,
\end{align*}
and%
\begin{gather}
\alpha_{h}:=\left\{
\begin{array}
[c]{ll}%
\sup\mathcal{A}_{h} & \text{if }\mathcal{A}_{h}\not =\emptyset,\\
a & \text{if }\mathcal{A}_{h}=\emptyset,
\end{array}
\right.  \quad\beta_{h}:=\left\{
\begin{array}
[c]{ll}%
\inf\mathcal{B}_{h} & \text{if }\mathcal{B}_{h}\not =\emptyset,\\
b & \text{if }\mathcal{B}_{h}=\emptyset,
\end{array}
\right. \nonumber\\
\underline{\theta}_{h}:=\min\left\{  \frac{1}{\beta_{h}-a},\frac{1}%
{b-\alpha_{h}}\right\}  ,\quad\overline{\theta}_{h}:=\frac{\alpha_{h}%
+\beta_{h}}{2}. \label{de}%
\end{gather}
Observe that, since $h\not \equiv 0$, $\underline{\theta}_{h}$ is well defined
and $\alpha_{h}<\beta_{h}$ (and so, $\overline{\theta}_{h}\in\left(
\alpha_{h},\beta_{h}\right)  $). Let us also set
\[
\delta_{\Omega}\left(  x\right)  :=dist\left(  x,\partial\Omega\right)
=\min\left(  x-a,b-x\right)  \text{.}%
\]

The next lemma provides some useful upper and lower bounds for $\mathcal{S}%
_{\phi}\left(  h\right)  $ when $h$ is nonnegative.

\begin{lemma}
\label{plap}Let $0\leq h\in L^{1}(\Omega)$ with $h\not \equiv 0$. Then in
$\overline{\Omega}$ it holds that%
\begin{equation}
\underline{\theta}_{h}\min\left\{  \int_{a}^{\overline{\theta}_{h}}\phi
^{-1}(\int_{y}^{\overline{\theta}_{h}}h)dy,\int_{\overline{\theta}_{h}}%
^{b}\phi^{-1}(\int_{\overline{\theta}_{h}}^{y}h)dy\right\}  \delta_{\Omega
}\leq\mathcal{S}_{\phi}\left(  h\right)  \leq\phi^{-1}(\int_{a}^{b}%
h)\delta_{\Omega}. \label{lap1}%
\end{equation}

\end{lemma}

\textit{Proof}. Let $v:=\mathcal{S}_{\phi}\left(  h\right)  $. Since
$\phi^{-1}$ is increasing and $h\geq0$ in $\Omega$, using (\ref{v1}) we see
that $v^{\prime}\left(  x\right)  =\phi^{-1}\left(  c_{h}-\int_{a}^{x}h\left(
t\right)  dt\right)  $ is nonincreasing and so $v$ is concave in $\Omega$.
Hence, since $v=0$ on $\partial\Omega$ and $v\not \equiv 0$ we deduce that
$v^{\prime}\left(  b\right)  <0<v^{\prime}\left(  a\right)  $ and therefore%
\begin{equation}
0<c_{h}<\int_{a}^{b}h\left(  t\right)  dt. \label{ch}%
\end{equation}
Employing again the fact that $\phi$ is increasing and (\ref{ch}) we find
that
\[
v^{\prime}\left(  a\right)  ,\left\vert v^{\prime}\left(  b\right)
\right\vert \leq\phi^{-1}\left(  \int_{a}^{b}h\right)
\]
and thus from the concavity of $v$ we obtain the second inequality in
(\ref{lap1}).

Let us prove the first inequality in (\ref{lap1}). We first claim that
\begin{equation}
v\geq\underline{\theta}_{h}\left\Vert v\right\Vert _{\infty}\delta_{\Omega
}\quad\text{in }\overline{\Omega}. \label{con}%
\end{equation}
In order to verify this, let $\xi\in\Omega$ be some point where $v$ reaches
its maximum (and so $v^{\prime}\left(  \xi\right)  =0$). We note that
$\xi>\alpha_{h}$. Indeed, when $\mathcal{A}_{h}=\emptyset$ this is obvious. If
$\mathcal{A}_{h}\not =\emptyset$, then by (\ref{v1}) we have $v\left(
x\right)  =\phi^{-1}\left(  c_{h}\right)  \left(  x-a\right)  $ for all
$x\in\left(  a,\alpha_{h}\right)  $, with $\phi^{-1}\left(  c_{h}\right)  >0$
by (\ref{ch}). In particular, $v$ is increasing for such $x$ and thus
$\xi>\alpha_{h}$ as asserted. Hence, recalling the concavity of $v$ we get
that for all $x\in\left[  \xi,b\right]  ,$%
\[
v\left(  x\right)  \geq\frac{v\left(  \xi\right)  \left(  b-x\right)  }{b-\xi
}\geq\frac{\left\Vert v\right\Vert _{\infty}}{b-\alpha_{h}}\delta_{\Omega
}\left(  x\right)  .
\]
Analogously, for $x\in\left[  a,\xi\right]  ,$%
\[
v\left(  x\right)  \geq\frac{v\left(  \xi\right)  \left(  x-a\right)  }{\xi
-a}\geq\frac{\left\Vert v\right\Vert _{\infty}}{\beta_{h}-a}\delta_{\Omega
}\left(  x\right)
\]
and the claim is proved.

Suppose now that $\xi\geq\overline{\theta}_{h}$. Taking into account that
$\phi$ is an homeomorphism with $\phi\left(  0\right)  =0$, that $v^{\prime
}\left(  x\right)  =\phi^{-1}\left(  c_{h}-\int_{a}^{x}h\right)  $ and
$v^{\prime}\left(  \xi\right)  =0$, we derive that $c_{h}=\int_{a}^{\xi}h$.
Then, recalling (\ref{v1}), that $\phi$ is increasing and $h\geq0$,
\begin{equation}
v\left(  \overline{\theta}_{h}\right)  =\int_{a}^{\overline{\theta}_{h}}%
\phi^{-1}\left(  \int_{a}^{\xi}h-\int_{a}^{y}h\right)  dy\geq\int
_{a}^{\overline{\theta}_{h}}\phi^{-1}\left(  \int_{y}^{\overline{\theta}_{h}%
}h\right)  dy. \label{vvh}%
\end{equation}
Assume now that $\xi\leq\overline{\theta}_{h}$. In this case we rewrite $v$
as
\[
v\left(  x\right)  =\int_{x}^{b}\phi^{-1}\left(  \widetilde{c}_{h}-\int
_{y}^{b}h\left(  t\right)  dt\right)  dy,
\]
where $\widetilde{c}_{h}$ is the unique constant such that $v\left(  a\right)
=0$. Moreover, reasoning as in the previous paragraph we see that
$\widetilde{c}_{h}=\int_{\xi}^{b}h$. Therefore,
\begin{equation}
v\left(  \overline{\theta}_{h}\right)  =\int_{\overline{\theta}_{h}}^{b}%
\phi^{-1}\left(  \int_{\xi}^{b}h-\int_{y}^{b}h\right)  dy\geq\int
_{\overline{\theta}_{h}}^{b}\phi^{-1}\left(  \int_{\overline{\theta}_{h}}%
^{y}h\right)  dy. \label{vvhh}%
\end{equation}
Taking into account (\ref{con}), (\ref{vvh}) and (\ref{vvhh}) we may infer the
first inequality in (\ref{lap1}) and this concludes the proof. $\blacksquare$

\bigskip

\begin{remark}
\label{hopf}Let $0\leq h\in L^{1}(\Omega)$ with $h\not \equiv 0$.\strut

\begin{enumerate}
\item[(i)] Observe that, since $\overline{\theta}_{h}\in\left(  \alpha
_{h},\beta_{h}\right)  $, the constant that appears in the first term of the
inequalities in (\ref{lap1}) is strictly positive.

\item[(ii)] For any $g\in\mathcal{C}(\overline{\Omega})$ with $g>0$ in
$\Omega$, note that $\alpha_{h}=\alpha_{hg}$ and $\beta_{h}=\beta_{hg}$.
Therefore, by the above lemma we have that%
\[
\mathcal{S}_{\phi}\left(  hg\right)  \geq\underline{\theta}_{h}\min\left\{
\int_{a}^{\overline{\theta}_{h}}\phi^{-1}(\int_{y}^{\overline{\theta}_{h}%
}hg)dy,\int_{\overline{\theta}_{h}}^{b}\phi^{-1}(\int_{\overline{\theta}_{h}%
}^{y}hg)dy\right\}  \delta_{\Omega}\quad\text{in }\overline{\Omega}.
\]

\end{enumerate}
\end{remark}

Let $f:\Omega\times\mathbb{R}\rightarrow\mathbb{R}$ be a Carath\'{e}odory
function (that is, $f\left(  \cdot,\xi\right)  $ is measurable for all $\xi
\in\mathbb{R}$ and $f\left(  x,\cdot\right)  $ is continuous for $a.e.$
$x\in\Omega$). Let $\mathcal{L}$ be as in (\ref{L}), and let us now consider
problems of the form%
\begin{equation}
\left\{
\begin{array}
[c]{ll}%
\mathcal{L}u=f\left(  x,u\right)  & \text{in }\Omega,\\
u=0 & \text{on }\partial\Omega.
\end{array}
\right.  \label{noli}%
\end{equation}
We say that $v\in\mathcal{C}(\overline{\Omega})$ is a \textit{subsolution }of
(\ref{noli}) if there exists a finite set $\Sigma\subset\Omega$ such that
$\phi(v^{\prime})\in AC_{loc}(\overline{\Omega}\,\backslash\,\Sigma),$
$v^{\prime}(\tau^{+}):=\lim_{x\rightarrow\tau^{+}}v^{\prime}(x)\in\mathbb{R}$,
$v^{\prime}(\tau^{-}):=\lim_{x\rightarrow\tau^{-}}v^{\prime}(x)\in\mathbb{R}$
for each $\tau\in\Sigma,$ and%
\begin{equation}
\left\{
\begin{array}
[c]{ll}%
\mathcal{L}v\leq f\left(  x,v\left(  x\right)  \right)  & a.e.\text{ }%
x\in\Omega,\\
v\leq0\text{ on }\partial\Omega, & v^{\prime}(\tau^{-})<v^{\prime}(\tau
^{+})\text{ for each }\tau\in\Sigma.
\end{array}
\right.  \label{subi}%
\end{equation}
If the inequalities in (\ref{subi}) are inverted, we say that $v$ is a
\textit{supersolution} of (\ref{noli}).

$\qquad$

For the reader's convenience we state the following existence theorem in the
presence of well-ordered sub and supersolutions (for a proof, see for instance
\cite[Theorem 7.16]{rach}).

\begin{theorem}
\label{subsup} Let $v$ and $w$ be sub and supersolutions respectively of
(\ref{noli}) such that $v\left(  x\right)  \leq w\left(  x\right)  $ for all
$x\in\overline{\Omega}$. Suppose there exists $g\in L^{1}\left(
\Omega\right)  $ such that
\[
\left\vert f\left(  x,\xi\right)  \right\vert \leq g\left(  x\right)
\quad\text{for }a.e.\text{ }x\in\Omega\text{ and all }\xi\in\left[  v\left(
x\right)  ,w\left(  x\right)  \right]  .
\]
Then there exists $u\in\mathcal{C}^{1}(\overline{\Omega})$ solution of
(\ref{noli}) with $v\leq u\leq w$ in $\overline{\Omega}$.
\end{theorem}

\section{Main results}

Before proving our main results, let us introduce the following conditions on
$\phi$ and $f$.\medskip

H1. There exist $t_{1}>0$ and an increasing homeomorphism $\psi$ defined in
$\left[  0,t_{1}\right]  $ such that $\psi\left(  0\right)  =0$ and%
\begin{equation}
\phi\left(  tx\right)  \leq\psi\left(  t\right)  \phi\left(  x\right)
\quad\text{for all }t\in\left[  0,t_{1}\right]  \text{, }x\geq0\text{.}
\label{h1}%
\end{equation}

H1'. There exists $p>0$ such that
\begin{gather}
\underline{\lim}_{t\rightarrow0^{+}}\frac{t^{p}}{\phi\left(  t\right)
}>0\text{,\quad and}\label{h32}\\
\overline{\lim}_{t\rightarrow\infty}\frac{\phi\left(  c_{\Omega}t\right)
}{\phi\left(  t\right)  }<\infty\text{,}\quad\text{where\quad}c_{\Omega
}:=\frac{b-a}{2}\text{.} \label{hbi}%
\end{gather}

H2. There exist $t_{2},M>0$ such that%
\begin{equation}
\phi\left(  tx\right)  \leq M\phi\left(  t\right)  \phi\left(  x\right)
\quad\text{for all }t\in\left[  0,t_{2}\right]  \text{, }x\in\left[
0,c_{\Omega}\right]  \text{.} \label{h31}%
\end{equation}

F1. There exist $\overline{t},k_{1},k_{2},q>0$ such that%
\begin{equation}
k_{1}t^{q}\leq f\left(  t\right)  \quad\text{for }t\in\left[  0,\overline
{t}\right]  \quad\text{and}\quad f\left(  t\right)  \leq k_{2}t^{q}%
\quad\text{for all }t\geq0. \label{fcsi}%
\end{equation}

F1'. There exist $\overline{t},k_{1},k_{2},q_{1},q_{2}>0$ such that
\begin{equation}
k_{1}t^{q_{1}}\leq f\left(  t\right)  \quad\text{for }t\in\left[
0,\overline{t}\right]  \quad\text{and}\quad f\left(  t\right)  \leq k_{2}%
\phi(t)^{q_{2}}\quad\text{for all }t\geq0. \label{f2}%
\end{equation}

We notice that $c_{\Omega}=\max_{\overline{\Omega}}\delta_{\Omega}$. Let us
also mention that the inequality in (\ref{h31}) appears (but for
\textit{large} values of $t$ and $x$) in the so-called $\Delta^{\prime}$
condition referred to Young functions (see e.g. \cite{rao}).

\begin{remark}
\label{tre}\strut

\begin{enumerate}
\item[(i)] Note that if $\left\vert \Omega\right\vert \leq2$ the condition
(\ref{hbi}) holds automatically since $\phi$ is increasing and thus in that
case H1' reduces to (\ref{h32}). On the other hand, if H1 is true with
$\psi\left(  t\right)  =ct^{p}$ for some $c,p>0$, fixing $x=1$ in (\ref{h1})
we see that H1 implies (\ref{h32}). In other words, in this particular case,
in \textquotedblleft small\textquotedblright\ domains H1 is stronger than H1'.
However, in general, these hypothesis are independent (see examples (a2) and
(d) at the end of the paper).

\item[(ii)] Suppose that $\phi$ fulfills H1' or H1 with $\psi\left(  t\right)
=ct^{p}$ for some $c,p>0$. Then the condition%
\begin{equation}
\underline{\lim}_{x\rightarrow0^{+}}\frac{\phi\left(  x\right)  }{x^{p}}>0
\label{as}%
\end{equation}
is sufficient in order for H2 to hold. Indeed, in any case we may assume
(\ref{h32}) (see (i)). Hence, given any $t_{0}>0$, there exists $M_{t_{0}}>0$
such that $\phi\left(  t\right)  \leq M_{t_{0}}t^{p}$ for all $t\in\left[
0,t_{0}\right]  $. Also, (\ref{as}) implies that for every $x_{0}>0$ there
exists $N_{x_{0}}>0$ such that $x^{p}\leq N_{x_{0}}\phi\left(  x\right)  $ for
all $x\in\left[  0,x_{0}\right]  $. It follows that for all $t\in\left[
0,1\right]  $ and $x\in\left[  0,c_{\Omega}\right]  $,%
\[
\phi\left(  tx\right)  \leq M_{c_{\Omega}}\left(  tx\right)  ^{p}\leq
M_{c_{\Omega}}N_{1}N_{c_{\Omega}}\phi\left(  t\right)  \phi\left(  x\right)
,
\]
and thus H2 is valid. We observe however that (\ref{as}) is not necessary for
H2 to be true (see examples (a4), (b) and (c) below).

\item[(iii)] Let us point out that if $\phi$ is differentiable in $\left(
0,c_{\Omega}\right)  $ and
\[
\sup_{t\in\left(  0,1\right)  ,\text{ }x\in\left(  0,c_{\Omega}\right)  }%
\frac{t\phi^{\prime}\left(  tx\right)  }{\phi\left(  t\right)  \phi^{\prime
}\left(  x\right)  }:=M<\infty\text{,}%
\]
then one can readily verify that H2 holds with $t_{2}=1$.

\item[(iv)] It is not difficult to check that the hypothesis H1 and H2 are
independent, and that the same is true for H1' and H2, see examples (a), (a2)
and (d).\bigskip
\end{enumerate}
\end{remark}

Our results shall provide us with solutions that lie in the interior of the
positive cone of $\mathcal{C}_{0}^{1}(\overline{\Omega}):=\{u\in
\mathcal{C}^{1}(\overline{\Omega}):u=0\text{ on }\partial\Omega\}$, which is
denoted by
\[
\mathcal{P}^{\circ}:=\left\{  u\in\mathcal{C}_{0}^{1}(\overline{\Omega
}):u>0\text{ in }\Omega\text{\ and }u^{\prime}\left(  b\right)  <0<u^{\prime
}\left(  a\right)  \right\}  .
\]

\begin{theorem}
\label{nuevo}Let $0\leq m\in L^{1}\left(  \Omega\right)  $ with $m\not \equiv
0$.

(i) Assume H1 and F1 with%
\begin{equation}
\underline{\lim}_{t\rightarrow0^{+}}\frac{t^{q}}{\psi\left(  t\right)
}=\infty\text{.} \label{nu}%
\end{equation}
Then for all $\lambda>0$ there exists $u=u_{\lambda}\in\mathcal{P}^{\circ}$
solution of (\ref{intro}).

(ii) Assume H1' and F1' with
\begin{equation}
q_{1}\in\left(  0,p\right)  \quad\text{and\quad}q_{2}\in\left(  0,1\right)
\text{.} \label{q12}%
\end{equation}
Then for all $\lambda>0$ there exists $u=u_{\lambda}\in\mathcal{P}^{\circ}$
solution of (\ref{intro}).

Moreover, in both (i) and (ii) it holds that
\begin{equation}
\lim_{\lambda\rightarrow0^{+}}\left\Vert u_{\lambda}\right\Vert _{\mathcal{C}%
^{1}(\overline{\Omega})}=0. \label{norma}%
\end{equation}

\end{theorem}

\begin{remark}
\label{lapla}When $\phi$ is the $p$-Laplacian, i.e. $\phi\left(  t\right)
=\left\vert t\right\vert ^{p-2}t$ with $p>1$, clearly H1 (with $\psi\left(
t\right)  =t^{p-1}$) and H1' (with $p-1$ in place of $p$ in (\ref{h32})) hold.
Furthermore, (\ref{nu}) is valid if and only if $q<p-1$, so in this case we
have the usual growth condition that characterizes the sublinear problems.
Observe also that, since for the $p$-Laplacian in (ii) we can take any
$q_{1}\in\left(  0,p-1\right)  $ and $1>q_{2}\approx1$, Theorem \ref{nuevo}
(i) and (ii) provide here the same result.
\end{remark}

\textit{Proof}. Let $\lambda>0$. We start proving (i). Let $\psi$,
$t_{1},\overline{t},k_{1},k_{2},q>0$ be given by H1 and F1 accordingly. By the
the continuity of $\phi^{-1}$ and the fact that $\phi^{-1}\left(  0\right)
=0$, there exists $\overline{\varepsilon}>0$ such that
\begin{equation}
\phi^{-1}(\varepsilon\int_{a}^{b}m\delta_{\Omega}^{q})\leq\frac{\overline{t}%
}{c_{\Omega}} \label{pri}%
\end{equation}
for all $\varepsilon\in\left(  0,\overline{\varepsilon}\right]  $, where
$c_{\Omega}$ is given by (\ref{hbi}). Also, let $\underline{\theta}_{m}$ and
$\overline{\theta}_{m}$ be as in (\ref{de}) and set%
\[
\mathcal{M}_{\Omega}:=\min\left\{  \int_{a}^{\overline{\theta}_{m}}\phi
^{-1}(\int_{y}^{\overline{\theta}_{m}}m\delta_{\Omega}^{q})dy,\int
_{\overline{\theta}_{m}}^{b}\phi^{-1}(\int_{\overline{\theta}_{m}}^{y}%
m\delta_{\Omega}^{q})dy\right\}  .
\]
It follows from the definition of $\overline{\theta}_{m}$ that $\mathcal{M}%
_{\Omega}>0$. Let us also write%
\[
M:=\max\left\{  \frac{1}{\lambda k_{1}\left(  \underline{\theta}%
_{m}\mathcal{M}_{\Omega}\right)  ^{q}},\lambda k_{2}(\phi^{-1}(\int_{a}%
^{b}m\delta_{\Omega}^{q}))^{q}\right\}  .
\]
We now observe that by (\ref{nu}) there exists $\varepsilon_{0}>0$ such that
$M\psi\left(  \varepsilon\right)  \leq\varepsilon^{q}$ for all $\varepsilon
\in\left[  0,\varepsilon_{0}\right]  $. Hence,
\begin{equation}
M\varepsilon\leq\psi^{-1}\left(  \varepsilon\right)  ^{q} \label{e0}%
\end{equation}
for $\varepsilon\in\left[  0,\psi\left(  \varepsilon_{0}\right)  \right]  $.
We notice next that H1 says that $t\phi^{-1}\left(  x\right)  \leq\phi
^{-1}\left(  \psi\left(  t\right)  x\right)  $ for all $t\in\left[
0,t_{1}\right]  $ and $x\geq0$, and therefore
\begin{equation}
\psi^{-1}\left(  r\right)  \phi^{-1}\left(  x\right)  \leq\phi^{-1}\left(
rx\right)  \label{y}%
\end{equation}
for all $r\in\left[  0,\psi\left(  t_{1}\right)  \right]  $ and $x\geq0$.

Let us choose
\begin{equation}
0<\varepsilon\leq\min\left\{  1,\overline{\varepsilon},\psi\left(
\varepsilon_{0}\right)  ,\psi\left(  t_{1}\right)  \right\}  \text{,}
\label{e2}%
\end{equation}
and for such $\varepsilon$ define $v:=\mathcal{S}_{\phi}\left(  \varepsilon
m\delta_{\Omega}^{q}\right)  $. Since $\varepsilon\leq\overline{\varepsilon}$
and $\delta_{\Omega}\leq c_{\Omega}$ in $\Omega$, the second inequality in
(\ref{lap1}) and (\ref{pri}) tell us that $\left\Vert v\right\Vert _{\infty
}\leq\overline{t}$. Consequently, taking into account (\ref{e0}), (\ref{y})
and (\ref{e2}), employing F1 and Remark \ref{hopf} (ii) we deduce that
\begin{gather}
\lambda m\left(  x\right)  f\left(  v\right)  \geq\lambda k_{1}m\left(
x\right)  v^{q}\geq\label{co2}\\
\lambda k_{1}m\left(  x\right)  \left[  \underline{\theta}_{m}\min\left\{
\int_{a}^{\overline{\theta}_{m}}\phi^{-1}(\varepsilon\int_{y}^{\overline
{\theta}_{m}}m\delta_{\Omega}^{q})dy,\int_{\overline{\theta}_{m}}^{b}\phi
^{-1}(\varepsilon\int_{\overline{\theta}_{m}}^{y}m\delta_{\Omega}%
^{q})dy\right\}  \delta_{\Omega}\right]  ^{q}\geq\nonumber\\
\lambda k_{1}m\left(  x\right)  \left(  \underline{\theta}_{m}\psi^{-1}\left(
\varepsilon\right)  \mathcal{M}_{\Omega}\delta_{\Omega}\right)  ^{q}%
\geq\varepsilon m\left(  x\right)  \delta_{\Omega}^{q}\left(  x\right)
=-\phi\left(  v^{\prime}\right)  ^{\prime}\text{\quad in }\Omega
\text{.}\nonumber
\end{gather}
In other words, $v$ is a subsolution of (\ref{intro}).

On the other side, we see that H1 yields that $\phi\left(  x\right)
/\psi\left(  t\right)  \leq\phi\left(  x/t\right)  $ for all $t\in\left(
0,t_{1}\right]  $ and $x\geq0$. Thus, $\phi^{-1}\left(  x/\psi\left(
t\right)  \right)  \leq\phi^{-1}\left(  x\right)  /t$ for such $t$ and $x$ and
so,
\begin{equation}
\phi^{-1}(\frac{x}{r})\leq\frac{\phi^{-1}\left(  x\right)  }{\psi^{-1}\left(
r\right)  } \label{z}%
\end{equation}
for all $r\in\left(  0,\psi\left(  t_{1}\right)  \right]  $ and $x\geq0$. Let
now $w:=\mathcal{S}_{\phi}\left(  \varepsilon^{-1}m\delta_{\Omega}^{q}\right)
$. Recalling (\ref{e0}), (\ref{e2}) and (\ref{z}) and utilizing again F1 and
Lemma \ref{plap}, we get that%
\begin{gather*}
\lambda m\left(  x\right)  f\left(  w\right)  \leq\lambda k_{2}m\left(
x\right)  w^{q}\leq\lambda k_{2}m\left(  x\right)  \left(  \phi^{-1}(\frac
{1}{\varepsilon}\int_{a}^{b}m\delta_{\Omega}^{q})\delta_{\Omega}\right)
^{q}\leq\\
\lambda k_{2}m\left(  x\right)  \left(  \frac{1}{\psi^{-1}\left(
\varepsilon\right)  }\phi^{-1}(\int_{a}^{b}m\delta_{\Omega}^{q})\delta
_{\Omega}\right)  ^{q}\leq\frac{1}{\varepsilon}m\left(  x\right)
\delta_{\Omega}^{q}\left(  x\right)  =-\phi\left(  w^{\prime}\right)
^{\prime}\text{\quad in }\Omega\text{,}%
\end{gather*}
and hence $w$ is a supersolution of (\ref{intro}). Moreover, since
$\varepsilon\leq1$ and $\mathcal{S}_{\phi}$ is nondecreasing (see Lemma
\ref{node}) we infer that $v\leq w$ in $\overline{\Omega}$. Then, we may apply
Theorem \ref{subsup} to obtain a solution $u_{\lambda}\in\mathcal{C}%
^{1}(\overline{\Omega})$ of (\ref{intro}) with $v\leq u_{\lambda}\leq w$ in
$\overline{\Omega}$, and since $v\in\mathcal{P}^{\circ}$ it also holds that
$u_{\lambda}\in\mathcal{P}^{\circ}$.

Let us prove (ii). Let $\overline{t},k_{1},k_{2},q_{1},q_{2}>0$ be given by
H1'. We note that (\ref{h32}) implies that $\phi\left(  t\right)  \leq Kt^{p}$
for all $t\in\left[  0,1\right]  $ and some $K>0$. Hence, we have that $t\leq
K\phi^{-1}\left(  t\right)  ^{p}$ for $t\in\left[  0,\phi\left(  1\right)
\right]  $, or equivalently,
\begin{equation}
\phi^{-1}\left(  t\right)  \geq\left(  t/K\right)  ^{1/p} \label{kk}%
\end{equation}
for such $t$. We now set%
\[
\mathcal{N}_{\Omega}:=\min\left\{  \int_{a}^{\overline{\theta}_{m}}(\int
_{y}^{\overline{\theta}_{m}}m\delta_{\Omega}^{q_{1}})^{1/p}dy,\int
_{\overline{\theta}_{m}}^{b}(\int_{\overline{\theta}_{m}}^{y}m\delta_{\Omega
}^{q_{1}})^{1/p}dy\right\}  >0,
\]
and similarly to (i) we define $v:=\mathcal{S}_{\phi}\left(  \varepsilon
m\delta_{\Omega}^{q_{1}}\right)  $, picking%
\begin{equation}
0<\varepsilon\leq\min\left\{  \overline{\varepsilon},\left(  \lambda
k_{1}\left(  \frac{\underline{\theta}_{m}\mathcal{N}_{\Omega}}{K^{1/p}%
}\right)  ^{q_{1}}\right)  ^{p/(p-q_{1})}\right\}  , \label{ew}%
\end{equation}
where $\overline{\varepsilon}>0$ is such that $\phi^{-1}(\overline
{\varepsilon}\int_{a}^{b}m\delta_{\Omega}^{q_{1}})\leq\overline{t}/c_{\Omega}%
$. As in the proof of (i) we have that $\left\Vert v\right\Vert _{\infty}%
\leq\overline{t}$. Thus, taking into account (\ref{q12}), (\ref{kk}) and
(\ref{ew}) and arguing as in (\ref{co2}) we derive that%
\begin{gather}
\lambda m\left(  x\right)  f\left(  v\right)  \geq\lambda k_{1}m\left(
x\right)  v^{q_{1}}\geq\label{co3}\\
\lambda k_{1}m\left(  x\right)  \left(  \underline{\theta}_{m}(\frac
{\varepsilon}{K})^{1/p}\mathcal{N}_{\Omega}\delta_{\Omega}\right)  ^{q_{1}%
}\geq\varepsilon m\left(  x\right)  \delta_{\Omega}^{q_{1}}\left(  x\right)
=-\phi\left(  v^{\prime}\right)  ^{\prime}\text{\quad in }\Omega
\text{.}\nonumber
\end{gather}

On the other hand, let $N:=\sup_{t>1}\phi\left(  c_{\Omega}t\right)
/\phi\left(  t\right)  <\infty$ (by (\ref{hbi})). For all $t\geq1$ we have
$\phi\left(  c_{\Omega}t\right)  \leq N\phi\left(  t\right)  $ and so%
\begin{equation}
\phi\left(  c_{\Omega}\phi^{-1}\left(  t\right)  \right)  \leq Nt \label{ome}%
\end{equation}
for all $t\geq\phi\left(  1\right)  $. Let $w:=\mathcal{S}_{\phi}\left(
\gamma m\right)  $ with%
\begin{equation}
\gamma\geq\max\left\{  \frac{\phi\left(  1\right)  }{\int_{a}^{b}m},\left(
\lambda k_{2}(N\int_{a}^{b}m)^{q_{2}}\right)  ^{1/(1-q_{2})}\right\}  .
\label{gaio}%
\end{equation}
Recalling F1', the upper bound given by Lemma \ref{plap} and that $q_{2}%
\in\left(  0,1\right)  $ and $\delta_{\Omega}\leq c_{\Omega}$ in $\Omega$,
employing (\ref{ome}) and (\ref{gaio}) we infer that
\begin{gather*}
\lambda m\left(  x\right)  f\left(  w\right)  \leq\lambda k_{2}m\left(
x\right)  \phi\left(  w\right)  ^{q_{2}}\leq\\
\lambda k_{2}m\left(  x\right)  \left(  \phi(\phi^{-1}(\gamma\int_{a}%
^{b}m)\delta_{\Omega})\right)  ^{q_{2}}\leq\lambda k_{2}m\left(  x\right)
\left(  \phi(c_{\Omega}\phi^{-1}(\gamma\int_{a}^{b}m))\right)  ^{q_{2}}\leq\\
\lambda k_{2}m\left(  x\right)  \left(  N\gamma\int_{a}^{b}m\right)  ^{q_{2}%
}\leq\gamma m\left(  x\right)  =-\phi\left(  w^{\prime}\right)  ^{\prime
}\text{\quad in }\Omega\text{.}%
\end{gather*}
Furthermore, enlarging $\gamma$ if necessary so that $\gamma\geq\varepsilon
c_{\Omega}^{q_{1}}$ and utilizing Lemma \ref{node} we can achieve that $w\geq
v$ in $\overline{\Omega}$ and thus we obtain a solution $u_{\lambda}%
\in\mathcal{P}^{\circ}$ of (\ref{intro}).

Finally, let us prove (\ref{norma}). Let $\lambda_{0}>0$ be fixed, and
consider $\lambda\in\left(  0,\lambda_{0}\right)  $. We first observe that the
solutions $u_{\lambda}$ obtained in either (i) or (ii) can be chosen such that
$\left\Vert u_{\lambda}\right\Vert _{\infty}\leq C$ with $C$ independent of
$\lambda$. Indeed, since $u_{\lambda_{0}}\in\mathcal{P}^{\circ}$ is a
supersolution of (\ref{intro}) for any $\lambda\in\left(  0,\lambda
_{0}\right)  $, and since the above part of the proof provides arbitrary small
subsolutions of (\ref{intro}) (that converge to $0$ in $\mathcal{C}%
(\overline{\Omega})$ as $\varepsilon\rightarrow0$, by the second inequality in
(\ref{lap1})), it follows from Theorem \ref{subsup} that there exist
$u_{\lambda}\in\mathcal{P}^{\circ}$ solutions of (\ref{intro}) such that
$0\leq u_{\lambda}\leq u_{\lambda_{0}}$ for all $\lambda\in\left(
0,\lambda_{0}\right)  $. So, $\left\Vert u_{\lambda}\right\Vert _{\infty}\leq
C$ as claimed. Taking into account this, the upper estimate in Lemma
\ref{plap} yields that%
\[
0\leq u_{\lambda}\left(  x\right)  =\mathcal{S}_{\phi}\left(  \lambda
mf\left(  u_{\lambda}\right)  \right)  \left(  x\right)  \leq\phi^{-1}\left(
\lambda\int_{a}^{b}mf\left(  u_{\lambda}\right)  \right)  \delta_{\Omega
}\left(  x\right)  \rightarrow0
\]
uniformly in $\overline{\Omega}$ as $\lambda\rightarrow0^{+}$ and so
$\lim_{\lambda\rightarrow0^{+}}\left\Vert u_{\lambda}\right\Vert _{\infty}=0$.

We choose next $\xi=\xi_{\lambda}\in\Omega$ such that $u_{\lambda}^{\prime
}\left(  \xi\right)  =0$. Integrating (\ref{intro}) over $\left(
a,\xi\right)  $ we get that $u_{\lambda}^{\prime}\left(  a\right)  =\phi
^{-1}\left(  \lambda\int_{a}^{\xi}mf\left(  u_{\lambda}\right)  \right)  $ and
hence by the above paragraph we see that $u_{\lambda}^{\prime}\left(
a\right)  \rightarrow0$ as $\lambda\rightarrow0^{+}$. Now, for any
$x\in\overline{\Omega}$, we integrate (\ref{intro}) over $\left(  a,x\right)
$ to find that
\[
u_{\lambda}^{\prime}\left(  x\right)  =\phi^{-1}\left(  \phi\left(
u_{\lambda}^{\prime}\left(  a\right)  \right)  +\lambda\int_{a}^{x}mf\left(
u_{\lambda}\right)  \right)  \rightarrow0
\]
uniformly when $\lambda\rightarrow0^{+}$. Thus, the proof of (\ref{norma}) is
complete. $\blacksquare$

\qquad

We next consider the case $r\in L^{1}\left(  \Omega\right)  $ with $r\geq0$,
that is, the problem%
\begin{equation}
\left\{
\begin{array}
[c]{ll}%
-\phi\left(  u^{\prime}\right)  ^{\prime}+r\left(  x\right)  \phi\left(
u\right)  =\lambda m\left(  x\right)  f\left(  u\right)  & \text{in }\Omega,\\
u=0 & \text{on }\partial\Omega.
\end{array}
\right.  \label{intror}%
\end{equation}

\begin{theorem}
\label{r}Let $0\leq m\in L^{1}\left(  \Omega\right)  $ with $m\not \equiv 0$.
Assume that $\phi$ fulfills H2, and suppose $\phi$ and $f$ satisfy the
hypothesis of Theorem \ref{nuevo} (i) or (ii), with $\psi\left(  t\right)
=ct^{p}$ for some $c,p>0$ in case (i). If either $r\leq m$ in $\Omega$ or
$m,r\in L^{\infty}\left(  \Omega\right)  $ and $\inf_{\Omega}m>0$, then for
all $\lambda>0$ there exists $u=u_{\lambda}\in\mathcal{P}^{\circ}$ solution of
(\ref{intror}). Moreover, these $u_{\lambda}$ satisfy (\ref{norma}).
\end{theorem}

\textit{Proof}. The proof follows the lines of the proof of Theorem
\ref{nuevo} and hence we only indicate the minor changes that are needed.

Let $\lambda>0$ and suppose the hypothesis of Theorem \ref{nuevo} (i) hold.
Let $t_{2},M>0$ be given by H2 and pick $\varepsilon>0$ such that
\[
\phi^{-1}(\varepsilon\int_{a}^{b}m\delta_{\Omega}^{q})\leq t_{2}.
\]
For such $\varepsilon$ define $v:=\mathcal{S}_{\phi}\left(  \varepsilon
m\delta_{\Omega}^{q}\right)  $. Taking $x=1$ in (\ref{h1}) (and recalling that
here $\psi\left(  t\right)  =ct^{p}$ for some $c,p>0$ ) we get that there
exists $K>0$ such that
\begin{equation}
\phi\left(  t\right)  \leq Kt^{p}\quad\text{for all }t\in\left[  0,c_{\Omega
}\right]  , \label{aau}%
\end{equation}
where $c_{\Omega}$ is given by (\ref{hbi}). Taking into account that
$\delta_{\Omega}\leq c_{\Omega}$ in $\Omega$, using Lemma \ref{plap} and
H2\ we derive that
\begin{equation}
\phi\left(  v\right)  \leq\phi\left(  \phi^{-1}(\varepsilon\int_{a}^{b}%
m\delta_{\Omega}^{q})\delta_{\Omega}\right)  \leq\varepsilon M\phi\left(
\delta_{\Omega}\right)  \int_{a}^{b}m\delta_{\Omega}^{q}\leq\varepsilon
MK\delta_{\Omega}^{p}\int_{a}^{b}m\delta_{\Omega}^{q}. \label{cofi}%
\end{equation}

Now, assume first that $r\leq m$ in $\Omega$. By (\ref{nu}) we have that
$q<p$. Thus, making $\varepsilon$ smaller if necessary, since $\psi
^{-1}\left(  t\right)  =\left(  t/c\right)  ^{1/p}$, from (\ref{co2}) and
(\ref{cofi}) we get that
\begin{gather*}
\lambda m\left(  x\right)  f\left(  v\right)  -r\left(  x\right)  \phi\left(
v\right)  \geq\\
m\left(  x\right)  \left(  \lambda k_{1}\left(  \underline{\theta}_{m}%
(\frac{\varepsilon}{c})^{1/p}\mathcal{M}_{\Omega}\delta_{\Omega}\right)
^{q}-\varepsilon MK\delta_{\Omega}^{p}\int_{a}^{b}m\delta_{\Omega}^{q}\right)
\geq\varepsilon m\left(  x\right)  \delta_{\Omega}^{q}\left(  x\right)
=-\phi\left(  v^{\prime}\right)  ^{\prime}.
\end{gather*}

On the other hand, if $r,m\in L^{\infty}\left(  \Omega\right)  $ and
$\underline{m}:=\inf_{\Omega}m>0$, for all $\varepsilon$ sufficiently small,
also from (\ref{co2}) and (\ref{cofi}) we deduce that
\begin{gather*}
\lambda m\left(  x\right)  f\left(  v\right)  -r\left(  x\right)  \phi\left(
v\right)  \geq\\
\lambda\underline{m}k_{1}\left(  \underline{\theta}_{m}(\frac{\varepsilon}%
{c})^{1/p}\mathcal{M}_{\Omega}\delta_{\Omega}\right)  ^{q}-\left\Vert
r\right\Vert _{\infty}\varepsilon MK\delta_{\Omega}^{p}\int_{a}^{b}%
m\delta_{\Omega}^{q}\geq\varepsilon\left\Vert m\right\Vert _{\infty}%
\delta_{\Omega}^{q}\left(  x\right)  \geq-\phi\left(  v^{\prime}\right)
^{\prime}.
\end{gather*}
Hence, in any case we obtain a subsolution of (\ref{intror}) which belongs to
$\mathcal{P}^{\circ}$. Furthermore, these subsolutions tend uniformly to zero
(by Lemma \ref{plap}) as $\varepsilon\rightarrow0$. Therefore, since the
solutions given by Theorem \ref{nuevo} (which also lie in $\mathcal{P}^{\circ
}$) are supersolutions of (\ref{intror}), Theorem \ref{subsup} yields the
desired solution $u_{\lambda}$. Moreover, it also follows that $\lim
_{\lambda\rightarrow0^{+}}\left\Vert u_{\lambda}\right\Vert _{L^{\infty
}(\Omega)}=0$, and similar computations to those in the last part of the proof
of Theorem \ref{nuevo} show that $u_{\lambda}$ satisfy (\ref{norma}).

Suppose now the assumptions of Theorem \ref{nuevo} (ii) hold. Then we set
$v:=\mathcal{S}_{\phi}\left(  \varepsilon m\delta_{\Omega}^{q_{1}}\right)  $,
where $q_{1}$ is given by H1'. Since (\ref{aau}) is true by (\ref{h32}),
proceeding as in (\ref{cofi}) we have that
\[
\phi\left(  v\right)  \leq\varepsilon MK\delta_{\Omega}^{p}\int_{a}^{b}%
m\delta_{\Omega}^{q_{1}}.
\]
Therefore, employing (\ref{co3}) in place of (\ref{co2}) and arguing as in the
above two paragraphs we can construct arbitrarily small subsolutions and thus
the proof can be completed as before. $\blacksquare$

\qquad

As a direct consequence of Theorems \ref{nuevo} and \ref{r} we are able to
provide an existence result also for
\begin{equation}
\left\{
\begin{array}
[c]{ll}%
\mathcal{L}u=\lambda m\left(  x\right)  f\left(  u\right)  & \text{in }%
\Omega,\\
u=0 & \text{on }\partial\Omega,
\end{array}
\right.  \label{corro}%
\end{equation}
where $\mathcal{L}u=-\phi\left(  u^{\prime}\right)  ^{\prime}$ or
$\mathcal{L}u=-\phi\left(  u^{\prime}\right)  ^{\prime}+r\left(  x\right)
\phi\left(  u\right)  $ accordingly, in the case where $m$ changes sign in
$\Omega$. As usual, we write $m=m^{+}-m^{-}$ with $m^{\pm}:=\max\left(  \pm
m,0\right)  $.

\begin{corollary}
\label{cori}Let $m\in L^{1}\left(  \Omega\right)  $ such that there exists an
open interval $\Omega_{0}\subset\Omega$ with $0\leq m\not \equiv 0$ in
$\Omega_{0}$. Suppose the hypothesis of one of the above theorems are
satisfied, with $m^{+}$ in place of $m$. Then for all $\lambda>0$ there exists
$u=u_{\lambda}\in\mathcal{C}^{1}(\overline{\Omega})$ nonnegative (and
nontrivial) solution of (\ref{corro}). Moreover, these $u_{\lambda}$ satisfy
(\ref{norma}).
\end{corollary}

\textit{Proof}. Let $\lambda>0$, and let $\overline{u}=\overline{u}_{\lambda
}\in\mathcal{P}^{\circ}$ be the solution of (\ref{corro}) with $m^{+}$ in
place of $m$, provided by some of the above theorems. It is clear that
$\overline{u}$ is a supersolution of (\ref{corro}).

On the other side, since $0\leq m\not \equiv 0$ in $\Omega_{0}$, an inspection
of the proofs of the aforementioned theorems show that we can find some
$z=z_{\lambda}\in\mathcal{C}^{1}(\overline{\Omega}_{0})$ with $z_{\lambda}%
\leq\overline{u}_{\lambda}$ in $\Omega_{0}$ and such that
\[
\left\{
\begin{array}
[c]{ll}%
\mathcal{L}z\leq\lambda m\left(  x\right)  f\left(  z\right)  & \text{in
}\Omega_{0},\\
z=0 & \text{on }\partial\Omega_{0}.
\end{array}
\right.
\]
Define now $\underline{u}_{\lambda}\in\mathcal{C}(\overline{\Omega})$ by
$\underline{u}_{\lambda}:=z_{\lambda}$ in $\Omega_{0}$ and $\underline
{u}_{\lambda}:=0$ in $\overline{\Omega}$\thinspace$\backslash\,\Omega_{0}$.
Then $\underline{u}_{\lambda}$ is a subsolution of (\ref{corro}) and this
yields the existence assertion.

To conclude the proof we note that the last assertion follows similarly to the
previous theorems, having in mind that $u_{\lambda}\leq\overline{u}_{\lambda}$
in $\Omega$ and that $\overline{u}_{\lambda}\rightarrow0$ uniformly as
$\lambda\rightarrow0^{+}$. $\blacksquare$

\qquad

\textbf{Examples}. We assume that $x\geq0$ since we may extend $\phi$
oddly.\smallskip

(a) Let $\varphi:\left[  0,\infty\right)  \rightarrow\left[  0,\infty\right)
$ be continuous and nondecreasing, with $\varphi$ increasing in $\left(
0,x_{0}\right)  $ for some $x_{0}>0$ if $\varphi\left(  0\right)  =0$. Define
\begin{equation}
\phi\left(  x\right)  :=x^{p}\varphi\left(  x\right)  ,\quad p>0. \label{ejj}%
\end{equation}
Then $\phi$ fulfills H1 with $t_{1}:=1$ and $\psi\left(  t\right)  :=t^{p}$
because
\[
\phi\left(  tx\right)  =\left(  tx\right)  ^{p}\varphi\left(  tx\right)
\leq\psi\left(  t\right)  x^{p}\varphi\left(  x\right)  =\psi\left(  t\right)
\phi\left(  x\right)
\]
for all $t\in\left[  0,1\right]  $ and $x\geq0$. Furthermore, in this case the
condition (\ref{nu}) of Theorem \ref{nuevo} is true if and only if $p>q$.

Let us note that here $\phi$ satisfies H2 if and only if $\varphi$ does. In
particular, taking some $\varphi$ which does \textit{not} fulfill H2 we obtain
a function $\phi$ that satisfies H1 but not H2 (one such $\varphi$ is for
instance $\varphi\left(  x\right)  =e^{-1/x}$ for $x>0$ and $\varphi\left(
0\right)  =0$).

We finally point out that if $\left\vert \Omega\right\vert \leq2$, the above
paragraph together with Remark \ref{tre} (i) imply the existence of some
$\phi$ which satisfies H1' but not H2.

Let us exhibit next some interesting particular cases:

(a1) Let
\[
\phi\left(  x\right)  :=x^{p_{1}}+x^{p_{2}},\quad p_{1}\geq p_{2}>0\text{.}%
\]
Since $\varphi\left(  x\right)  :=\phi\left(  x\right)  /x^{p_{2}}$ is
increasing, by the first paragraph in (a) we get that H1 holds, and it is also
clear that H2 is true with $M=1$ and any $t_{2}>0$.

(a2) Let%
\[
\phi\left(  x\right)  :=e^{x^{p}}-1,\quad p>0.
\]
A brief computation shows that $\phi\left(  x\right)  /x^{p}$ is increasing
and so $\phi$ fulfills H1. Moreover, taking into account Remark \ref{tre} (ii)
we see that $\phi$ satisfies H2. Note also that if $\left\vert \Omega
\right\vert >2$, then (\ref{hbi}) is not valid. In particular, $\phi$
\textit{does not} fulfill H1' in this case (and for any $\Omega$, $\phi$
\textit{neither }satisfies the conditions $\left(  \Phi\right)  $ and $\left(
\Phi^{\prime}\right)  $ at the introduction nor the one in \cite{xu}).

(a3) Let
\[
\phi\left(  x\right)  :=e^{x}-x-1.
\]
One can check that $\phi\left(  x\right)  /x^{2}$ is increasing and therefore
H1 holds. Also, recalling Remark \ref{tre} (ii) we deduce that $\phi$
satisfies H2. We observe that $\phi$ \textit{does not }satisfy the assumptions
$\left(  \Phi\right)  $ and $\left(  \Phi^{\prime}\right)  $ or the one in
\cite{xu}.

(a4) Let%
\[
\phi\left(  x\right)  :=\frac{x^{p_{1}}}{1+x^{p_{2}}},\quad p_{1}>p_{2}>0.
\]
Since $\phi\left(  x\right)  /x^{p_{1}-p_{2}}$ is increasing we infer that H1
is valid. Although (\ref{as}) \textit{does not }hold with $p=p_{1}-p_{2}$, one
can verify that $\phi$ fulfills H2 with $t_{2}:=1$ and $M:=2\left(
1+c_{\Omega}^{p_{2}}\right)  $.

(b) Let
\[
\phi\left(  x\right)  :=x\left(  \left\vert \ln x\right\vert +1\right)  .
\]
It can be proved that $\phi$ \textit{cannot} be written as in (\ref{ejj}) with
$p>0$ and $\varphi$ nondecreasing. Let us demonstrate, however, that $\phi$
satisfies H1 and H2. Let $p\in\left(  0,1\right)  $. We choose $t_{1}>0$ such
that $\left\vert \ln t\right\vert \leq1/t^{1-p}-1$ for all $t\in\left[
0,t_{1}\right]  $. Then, for such $t$ and all $x\geq0$ we have that%
\begin{gather}
\phi\left(  tx\right)  \leq tx\left(  \left\vert \ln t\right\vert +\left\vert
\ln x\right\vert +1\right)  \leq\label{tx}\\
tx\left[  \left(  1/t^{1-p}-1\right)  \left(  \left\vert \ln x\right\vert
+1\right)  +\left\vert \ln x\right\vert +1\right]  =t^{p}\phi\left(  x\right)
\nonumber
\end{gather}
and H1 holds. Also, employing the first inequality in (\ref{tx}) it is easy to
see that H2 is true with $M=1$ and any $t_{2}>0$. Let us finally note that
(\ref{as}) \textit{is not }valid with any $p\in\left(  0,1\right)  $.

(c) Let%
\[
\phi\left(  x\right)  :=x-\ln\left(  x+1\right)  .
\]
Then $\phi$ fulfills (\ref{h32}) with $p=1$ and also (\ref{hbi}). In other
words, H1' holds (let us remark that, despite it is less direct, one can prove
that $\phi$ also satisfies H1 with $t_{1}=1$ and $\psi\left(  t\right)  =ct$
for some $c>0$ large enough). Although (\ref{as}) \textit{does not }hold with
$p=1$, from Remark \ref{tre} (iii) we deduce that H2 is valid since%
\[
\frac{t\phi^{\prime}\left(  tx\right)  }{\phi\left(  t\right)  \phi^{\prime
}\left(  x\right)  }=\frac{t^{2}\left(  x+1\right)  }{\left(  t-\ln\left(
t+1\right)  \right)  \left(  tx+1\right)  }\leq\frac{t^{2}\left(  c_{\Omega
}+1\right)  }{\left(  t-\ln\left(  t+1\right)  \right)  }\leq\frac{c_{\Omega
}+1}{1-\ln2}%
\]
for all $t\in\left(  0,1\right)  $ and all $x\in\left(  0,c_{\Omega}\right)  $.

(d) Let
\[
\phi\left(  x\right)  :=\left(  \ln\left(  x+1\right)  \right)  ^{p},\quad
p>0.
\]
One can readily check (\ref{h32}) and (\ref{hbi}), and thus H1' is true. Also,
utilizing again Remark \ref{tre} (ii) we get that H2 holds. Let us observe
that $\phi$ \textit{does not }satisfy H1 (and hence \textit{neither }fulfills
the conditions $\left(  \Phi\right)  $ or $\left(  \Phi^{\prime}\right)  $ at
the introduction) because
\[
\lim_{x\rightarrow\infty}\frac{\phi\left(  tx\right)  }{\phi\left(  x\right)
}=1\quad\text{for all }t>0
\]
and so there is not a continuous $\psi$ such that (\ref{h1}) is valid and
$\psi\left(  0\right)  =0$. Furthermore, this tell us that $\phi$ neither
meets the assumptions in \cite{xu}.

\end{document}